
\baselineskip=14pt
\parskip=10pt

\font\eightrm=cmr8 

\magnification=\magstephalf

\def\1{{\overline{1}}}
\def\2{{\overline{2}}}
\parindent=0pt
\overfullrule=0in

\def\frac#1#2{{#1 \over #2}}
\centerline
{\bf 
Experimenting with the Dym-Luks Ball and Cell Game (almost) Sixty Years Later 
}
\bigskip
\centerline
{\it Shalosh B. EKHAD and Doron ZEILBERGER}
\centerline
\qquad \qquad \qquad 
{\it Dedicated to Harry Dym (b. Jan. 26, 1938) on  his  eighty-fifth birthday}
\bigskip

{\bf Abstract}:  
This is a symbolic-computational redux, and extension, of a beautiful paper, by Harry Dym and Eugene Luks, published in 1966 (but written in 1964) about
a certain game with balls and cells.

{\bf Preface}

When Harry Dym and Eugene Luks were graduate students at MIT (working with Henry P. McKean, Jr. and Kenkichi Iwasawa, respectively) they
collaborated on a ``{\it fun}'' paper [DL] not directly related to their dissertation topics. Here is how they introduced it.

{\it
``Each of $r$ balls is placed at random into one of $n$ cells. A ball is considered
``captured" if (after all r balls have been distributed) it is the sole occupant of its
 cell. Captured balls are eliminated from further play. This completes the first
 ``trial." The remaining balls are recovered and the process repeated (trials,
 $2, 3, 4, \dots $ etc.). The play continues until all balls have been captured. The
 number of trials required to achieve this state is called the duration of the game.''
}

In this modest tribute, we implement (in Maple) and {\it extend} their beautiful paper.

This will show the power of {\it symbolic computation}. We will  describe a Maple package, {\tt DymLuks.txt}, written by DZ, and diligently executed by SBE,
available from the front of this article:

{\tt https://sites.math.rutgers.edu/\~{}zeilberg/mamarim/mamarimhtml/dym.html} \quad .

This web-page also has   numerous input and output files referred to in this article.

Dym and Luks only considered the
{\it expected duration}, but we will go further and compute the  {\it variance},  and higher moments. We also compute
the full {\it probability generating function} for {\it any} specific number of balls and number of cells.
Also for a {\it fixed} number of balls $r$, we find explicit expressions, in $n$, for these quantities. 

{\bf A Quick summary of the Dym-Luks paper}

Dym and Luks viewed the game as a {\bf Markov process}. Fix $n$ (the number of cells) once and for all. If currently there are $r$ balls, then after one iteration,
assuming that $t$ balls were removed, there are  $r-t$ balls  that remain, where $0 \leq t \leq r$. 
Calling this probability $P_{r,r-t}(n)$,  invoking a clever inclusion-exclusion argument, they derived ([DL], p. 517)
$$
P_{r,r-t}(n) \, = \,
\sum_{j=t}^{n} (-1)^{j-t} {{j} \choose {t}} \, {{n} \choose {j}}\, {{r} \choose {j}}\, j! \frac{(n-j)^{r-j}}{n^r} \quad .
$$

[{\eightrm This is implemented in procedure {\tt Prt(n,r,t)} in our Maple package DymLuks.txt. Note that {\tt Prt(n,r,t)} is $P_{r,r-t}(n)$.}]

Then they focused on the case of a fixed number of cells, $n$, and arbitrary number of balls $r$, and looked at the behavior of $M_n(r)$, the expected duration, and proved that, for any fixed $n$,
as $r \rightarrow \infty$, we have
$$
M_n(r) \, = \, \sum_{j=1}^{r} j^{-1} \, \left( \frac{n}{n-1} \right)^{j-1} \, + \, O(1) \quad .
$$

{\bf Probability Generating Functions and Moments for the Dym-Luks Ball and Cell Game}

Rather than just talk about the expectation, 
we will compute the full {\bf probability generating function}, let's call it, $F_{r,n}(x)$.
This is the rational function  whose coefficient of $x^i$ (in its Maclaurin expansion) is
the {\bf exact} probability that the game will terminate in {\bf exactly} $i$ rounds. 
Once we have it, we easily get the expectation, $M_n(r)$, that equals $F_{n,r}'(1)$, and
the higher moments (see below).
We have 
$$
F_{r,n}(x) \, = \,x \left ( \sum_{t=0}^{r} P_{r,r-t} F_{r-t,n}(x) \right ) \quad .
$$
Hence we have the {\bf recurrence}
$$
F_{r,n}(x) \, = \,\frac{x}{1-P_{r,r}x} \left ( \sum_{t=1}^{r} P_{r,r-t} F_{r-t,n}(x) \right ) \quad .
$$

[{\eightrm This is implemented in procedure {\tt GFrn(r,n,x)} in our DymLuks.txt}]. For each {\bf numeric} $r$ and $n$ we get a certain
rational function in $x$. For example, 
$$
F_{1,1}(x) \, = \, x \quad, \quad
F_{2,2}(x) \, = \, -\frac{x}{-2+x} \quad, \quad
F_{3,3}(x) \, = \, \frac{2 x \left(5 x +3\right)}{\left(-3+x \right) \left(-9+x \right)}
\quad, \quad
$$
$$
F_{4,4}(x) \, = \, -\frac{3 x \left(25 x^{2}+316 x +64\right)}{\left(-4+x \right) \left(-16+x \right) \left(-32+5 x \right)} \quad,
$$
$$
F_{5,5}(x) \, = \, \frac{24 x \left(767 x^{3}+63115 x^{2}+182125 x +15625\right)}{\left(-5+x \right) \left(-25+x \right) \left(-125+13 x \right) \left(-625+41 x \right)} \quad .
$$

To see the expressions for $F_{r,r}(x)$ for $r \leq 40$ (always rational functions of $x$, of course), see the output file:

{\tt https://sites.math.rutgers.edu/\~{}zeilberg/tokhniot/oDymLuks3.txt} \quad .

From these we can deduce the diagonal sequence $\{M_r(r)\}$. It is not clear to us, with the available data,  whether this sequence tends to a `universal'  constant, or whether it (very!) slowly increases to infinity.
At any rate, we are almost sure that there is a {\it limiting distribution} (once you scale it). What is it?

To get an idea how this sequence starts, look at the output file

{\tt https://sites.math.rutgers.edu/\~{}zeilberg/tokhniot/oDymLuks5.txt} \quad ,

that also gives the sequence of variances.

{\bf Probability Generating Functions for the Duration with a Fixed number of balls and General Number of Cells}

When $r$ is fixed , and we let $n$ vary, we can get closed-form expressions, as rational functions in $x$ {\bf and} $n$,for the general
$F_{r,n}(x)$. For example

$$
F_{1,n}(x) \, = \, x \quad, \quad
F_{2,n}(x) \, = \, \frac{x \left(n -1\right)}{n -x} \quad, \quad
F_{3,n}(x) \, = \, 
\frac{x \left(n^{3}+2 n^{2} x -3 n^{2}-3 n x +2 n +x \right)}{\left(n -x \right) \left(n^{2}-x \right)} \quad,
$$
$$
F_{4,n}(x) \, = \, 
\frac{x \left(n^{6}+5 n^{5} x -6 n^{5}-15 n^{4} x +3 n^{3} x^{2}+11 n^{4}+9 n^{3} x -2 n^{2} x^{2}-6 n^{3}+3 n^{2} x -3 n \,x^{2}-2 n x +2 x^{2}\right)}{\left(n -x \right) \left(n^{2}-x \right) \left(n^{3}-3 n x +2 x \right)} \quad ,
$$
$$
F_{5,n}(x) \, = \, 
$$
$$
\frac{x}{\left(n -x \right) \left(n^{2}-x \right) \left(n^{3}-3 n x +2 x \right) \left(n^{4}-10 n x +9 x \right)} \cdot
$$
$$
(n^{10}+9 n^{9} x -10 n^{9}-39 n^{8} x +59 n^{7} x^{2}+35 n^{8}-13 n^{7} x -284 n^{6} x^{2}+12 n^{5} x^{3}
$$
$$
-50 n^{7}+280 n^{6} x +508 n^{5} x^{2}-38 n^{4} x^{3}+24 n^{6}-494 n^{5} x 
$$
$$
-415 n^{4} x^{2}+35 n^{3} x^{3}+359 n^{4} x +111 n^{3} x^{2}+20 n^{2} x^{3}-102 n^{3} x +39 n^{2} x^{2}-47 n \,x^{3}-18 n \,x^{2}+18 x^{3})  \quad .
$$
For expressions for $F_{r,n}(x)$ for $r \leq 40$, please consult the (large!) output file

{\tt https://sites.math.rutgers.edu/\~{}zeilberg/tokhniot/oDymLuks3.txt } \quad .

{\bf Moments}

Once we have an explicit expression for the probability generating function (always a rational function) we can let the
computer compute, for each specific $r$, but general $n$, not only the expectation (the only concern in [DL]), but also the
variance, and higher moments. Recall that the average, called $M_n(r)$ in [DL] is
$$
M_n(r)= x \frac{d}{dx} F_{r,n}(x)|_{x=1} = F_{r,n}'(1) \quad .
$$
Calling our random variable $X_{n,r}$ (so $M_n(r)=E[X_{n,r}]$), we have
$$
E[X_{n,r}^i] \, = \,
(x \frac{d}{dx})^i F_{r,n}(x)|_{x=1}  \quad ,
$$
and from this the computer can easily find the {\it moments about the mean}.
Recall that  the $i^{th}$ moment about the mean, $m_i$, is  $E[(X_{n,r}-M_n(r))^i]$. In particular the second moment about the mean is the {\it variance}, $Var_n(r)$.
Also recall that the {\it scaled} moments (about the mean)  are $m_i/m_2^{i/2}$.

Here are the first few expressions for $M_n(r)$ for small $r$.
$$
M_n(1) \, = \, 1 \quad,
$$
of course, followed by:
$$
M_n(2) \, = \, \frac{n}{n -1} \quad, \quad
M_n(3) \, = \,\frac{n \left(n +3\right)}{n^{2}-1} \quad, \quad
M_n(4) \, = \, 
\frac{\left(n^{2}+7 n -2\right) n^{2}}{\left(n^{3}-3 n +2\right) \left(n +1\right)} \quad,
$$
$$
M_n(5) \, = \,
\frac{n^{2} \left(n^{5}+12 n^{4}-6 n^{3}+48 n^{2}-125 n +10\right)}{\left(n^{4}-10 n +9\right) \left(n^{2}+n -2\right) \left(n +1\right)} \quad.
$$
Not surprisingly they all tend to $1$ as $n$ goes to infinity. After all, if you have many cells and only a few balls, they are all likely to lend in different cells,
making them all removable.

For the explicit expressions for all $r\leq 20$, as well as expressions for the variance, {\it skewness}, and {\it kurtosis} (i.e. the scaled third and fourth moments, respectively), see the output file

{\tt https://sites.math.rutgers.edu/\~{}zeilberg/tokhniot/oDymLuks1t.txt} \quad .

{\bf Fixed number of Cells and General Number of Balls}

Except when there are two cells (i.e. $n=2$), there is no {\it closed form} expression for $F_{r,n}(x)$ and even not for $M_r(n)$, but following Dym and Luks we can get
very good approximations. It is easy to see that, with a fixed $n$, we have
$$
P_{r,r}(n) = 1 - r (\frac{n-1}{n})^{r-1}  \cdot (1+ O(\beta_1^r)) \quad,
$$
$$
P_{r,r-1}(n) = r (\frac{n-1}{n})^{r-1} \cdot  (1+ O(\beta_2^r)) \quad,
$$
$$
P_{r,r-i}(n)= r (\frac{n-1}{n})^{r-1} \, \cdot O(\beta_3^r) \quad, \quad i \geq 2 \quad,
$$
where  $0<\beta_1,\beta_2,\beta_3<1$. In other words, we can approximate the process by looking at a very simplified Markov process where the 
``particle" either moves one unit down (i.e. from $r$ to $r-1$) with probability  $r (\frac{n-1}{n})^{r-1}$,
and stays in place otherwise. Let's consider the more general situation where for an {\it arbitrary} sequence $a(r)$,
the particle either stays where it currently is, position $r$ say, or goes down to $r-1$ with probability $a(r)$.
For the Ball and Cell Game, and fixed number of cells $n$, we have $a(r)=r\,(\frac{n-1}{n})^{r-1}$.

Let's consider this more general scenario. Having fixed the sequence $a(r)$, let $F_r(x)$ be the probability generating function for the duration.
We have:
$$
F_r(x)= x \cdot \left( (1-a(r))\, F_r(x) + a(r) F_{r-1}(x) \right) \quad,
$$
hence
$$
F_r(x) \, = \, \frac{a(r)\,x}{1-x(1-a(r))} \cdot F_{r-1}(x) \quad ,
$$
that, in turn, implies that
$$
F_r(x) \, = \, \prod_{i=1}^{r} \frac{a(i)\,x}{1-x(1-a(i))}  \quad .
$$
By taking derivative, and then plugging-in $x=1$ we  immediately get that the expectation is
$$
\sum_{i=1}^{r} \frac{1}{a(i)} \quad .
$$
By taking the second derivative and doing some simple manipulations, we get that the variance is
$$
\sum_{i=1}^{r} \frac{1}{a(i)^2} - \sum_{i=1}^{r} \frac{1}{a(i)} \quad .
$$
We can get similar expressions for the higher moments, but for general $a(r)$ they won't be very useful.
In the Dym-Luks case there is no ``closed-form", but if we consider the closely analogous case of
$a(r)=\alpha^r$, with $0<\alpha<1$, one (or rather one's computer) gets explicit expressions, not only for the expectation, but also for the variance and higher moments.

{\bf Explicit Expressions for the Moments of the Longevity of the Markov Process Where $r$ goes to $r-1$ with probability $\alpha^r$ and stays at $r$ with probability $1-\alpha^r$}

We have that the expectation is
$$
\frac{1-\alpha^r}{(1-\alpha)\alpha^r} \quad,
$$
and the variance is
$$
\frac{(1-\alpha^r)\,(1-\alpha^{r+1})}{(1-\alpha^2)\,\alpha^{2r}} \quad .
$$

In particular, the {\it coefficient of variation} (the  standard deviation divided by the expectation) converges to
$$
\sqrt{\frac{1-\alpha}{1+\alpha}} \quad .
$$

Explicit formulas for the general higher moments, up to the tenth, can be gotten from:

{\tt https://sites.math.rutgers.edu/\~{}zeilberg/tokhniot/oDymLuks7.txt} \quad .

Let's just state the limiting expressions as $r$ goes to infinity.

The limiting {\it skewness} (i.e. the scaled third-moment about the mean), as $r$ goes to infinity, is
$$
\sqrt{\frac{4 \left(1- \alpha\right) \left(\alpha +1\right)^{3}}{\left(\alpha^{2}+\alpha +1\right)^{2}}} \quad .
$$

The limiting {\it kurtosis} (i.e. the scaled fourth-moment about the mean) is
$$
\frac{3 \left(3-\alpha^{2}\right)}{\alpha^{2}+1}
\quad .
$$

The limiting {\it scaled fifth-moment}  (about the mean) , as $r$ goes to infinity, is
$$
\sqrt{\frac{16 \left(1-\alpha \right) \left(\alpha^{4}+\alpha^{3}-5 \alpha^{2}-11 \alpha -11\right)^{2} \left(\alpha +1\right)^{3}}{\left(\alpha^{2}+\alpha +1\right)^{2} \left(\alpha^{4}+\alpha^{3}+\alpha^{2}+\alpha +1\right)^{2}}}
\quad .
$$

The limiting {\it scaled sixth-moment}  (about the mean) , as $r$ goes to infinity, is
$$
\frac{5 \alpha^{8}+5 \alpha^{7}-45 \alpha^{6}-130 \alpha^{5}-180 \alpha^{4}-50 \alpha^{3}+135 \alpha^{2}+265 \alpha +265}{\left(\alpha^{2}+1\right) \left(\alpha^{2}-\alpha +1\right) \left(\alpha^{2}+\alpha +1\right)^{2}} \quad .
$$

{\bf How good are the Luks-Dym Approximations of $M_n(r)$?}

At the very end of [DL], the authors define the error
$$
E_n(r) := M_n(r) \, \, - \, \, \sum_{j=1}^{r} j^{-1} \, \left ( \frac{n}{n-1} \right )^{j-1}  \quad.
$$
Of course, as they noted, $E_2(r)=0$, and with a 1964 computer, they found out that $|E_3(r)| \leq 0.25$, and that it approaches $0.042$ as $r$ approaches infinity. Using simulation they
noticed that $E_n(r)$ seems to be less than one in magnitude, at least in the range $r \leq 5n$. With a 2023 computer we corroborated this. See the output file.

{\tt https://sites.math.rutgers.edu/\~{}zeilberg/tokhniot/oDymLuks4.txt} \quad .

We estimate that
$$
\lim_{r \rightarrow \infty} E_3(r) \approx .04213658385 \quad,
$$
(as already computed in [DL])
$$
\lim_{r \rightarrow \infty} E_4(r) \approx  .254461 \quad,
$$
$$
\lim_{r \rightarrow \infty} E_5(r) \approx  .5312 \quad \quad .
$$

Since $M_n(r)$ is so large, a more meaningful measure of the quality of the approximation is not the difference between $M_n(r)$ and its approximation, $\sum_{j=1}^{r} \, j^{-1}  \left (  \frac{n}{n-1} \right )^{j-1}$,
but rather the ratio. See the output file

{\tt https://sites.math.rutgers.edu/\~{}zeilberg/tokhniot/oDymLuks4R.txt} \quad ,

that also contains ratios of the variance with its approximation
$\sum_{j=1}^{r} j^{-2} \, (\frac{n}{n-1} )^{2j-2}-\sum_{j=1}^{r} j^{-1} \, \left ( \frac{n}{n-1} \right )^{j-1}$. As you can see, they all converge to $1$ very fast.

{\bf Simulation}

Out Maple package also has simulation procedures. Try {\tt DL(r,n)}, and {\tt DLv(r,n)} for a verbose version. To see ten examples starting with $1000$ balls and
$1000$ cells, consult the output file

{\tt https://sites.math.rutgers.edu/\~{}zeilberg/tokhniot/oDymLuks6.txt} \quad .

{\bf Conclusion}

Happy  (early) $85^{th}$ birthday, Harry, and Happy (late) $83^{th}$ birthday, Gene. Keep up the good work!

{\bf References}

[DL] Harry Dym and Eugene M. Luks, {\it  On the mean duration of a Ball and Cell game; a first  passage problem},
The Annals of Mathematical Statistics {\bf 37} (1966),  517-521. \hfill\break
{\tt https://sites.math.rutgers.edu/\~{}zeilberg/akherim/DymLuks.pdf} \quad .

\bigskip
\hrule
\bigskip
Shalosh B. Ekhad and Doron Zeilberger, Department of Mathematics, Rutgers University (New Brunswick), Hill Center-Busch Campus, 110 Frelinghuysen
Rd., Piscataway, NJ 08854-8019, USA. \hfill\break
Email: {\tt  ShaloshBEkhad at gmail dot com} \quad, \quad {\tt DoronZeil] at gmail dot com}   \quad .

Written: {\bf Jan. 12, 2023}.

\end